\documentclass[letterpaper, 10 pt, conference]{ieeeconf}

\usepackage{cite}
\usepackage{graphicx}
\usepackage{amsfonts}
\usepackage{amsmath} 
\usepackage{amssymb}
\usepackage{enumerate}
\usepackage{mathrsfs}
\usepackage{cleveref}
\usepackage[labelformat=simple]{subcaption}

\usepackage{epstopdf}
\usepackage{tikz}
\usetikzlibrary{arrows}

\usepackage{todonotes}
\newcommand{\hnorm}[1]{\left\|#1\right\|}
\usepackage[normalem]{ulem}

\graphicspath{{figures/}}

\usepackage{geometry}
\newtheorem{thm}{Theorem}

\newtheorem{cor}{Corollary}
\newtheorem{defn}{Definition}

\newtheorem{ass}{Assumption}

\newcommand{\R}{{\mathbb R}}

\newcommand{\prob}[1]{\mathop{{\rm Pr}} \left(#1 \right)}

\newcommand{\Nset}{\mathbb{N}}

\newcommand{\Rset}{\mathbb{R}}

\newcommand{\Zset}{\mathbb{Z}}

\newcommand{\Ucal}{{\cal U}}

\newcounter{l1}
\newcounter{l2}
\newcounter{l3}
\newcommand{\bdotlist}{\begin{list}{$\bullet$}{}}
\newcommand{\bboxlist}{\begin{list}{$\Box$}{}}
\newcommand{\bbboxlist}{\begin{list}{\raisebox{.005in}{{\tiny
$\blacksquare$ \ \ }}}{}}
\newcommand{\bdashlist}{\begin{list}{$-$}{} }
\newcommand{\blist}{\begin{list}{}{} }
\newcommand{\barablist}{\begin{list}{\arabic{l1}}{\usecounter{l1}}}
\newcommand{\balphlist}{\begin{list}{(\alph{l2})}{\usecounter{l2}}}
\newcommand{\bAlphlist}{\begin{list}{\Alph{l2}.}{\usecounter{l2}}}
\newcommand{\bdiamlist}{\begin{list}{$\diamond$}{}}
\newcommand{\bromalist}{\begin{list}{(\roman{l3})}{\usecounter{l3}}}

\newcommand{\rand}{\boldsymbol{u}}
\newcommand{\bzeta}{\boldsymbol{\zeta}}

\newcommand{\risk}{1-\alpha}
\newcommand{\conf}{\alpha}

\newcommand{\lev}[1]{\mathcal{U}_{\phi}(#1)}

\newcommand{\uncPar}{\rand}
\newcommand{\uncRel}{u}
\newcommand{\relLev}{\conf}
\newcommand{\vioLev}{\risk}
\newcommand{\uncSet}{\mathcal{U}}
\newcommand{\batNum}{n}

\newcommand{\batSiz}{m}
\newcommand{\batSizMem}{i}

\newcommand{\buf}{\varepsilon}
\newcommand{\confLev}{\delta}

\IEEEoverridecommandlockouts

\begin{document}

\title{\LARGE \bf Nonparametric Estimation of Uncertainty Sets for Robust Optimization}

\author{Polina Alexeenko$^{\dagger}$ and Eilyan Bitar$^{\dagger}$
\thanks{This work was supported in part by NSF grant ECCS-1351621, NSF grant IIP-1632124, and an NSF Graduate Research Fellowship.}
\thanks{$\dagger$ Polina Alexeenko  and Eilyan Bitar  are with the School of Electrical and Computer Engineering, Cornell Univ., Ithaca, NY, 14853, USA.  Emails: {\tt\small \{pa357, eyb5\}@cornell.edu}}%
}

\maketitle 

\begin{abstract}
We investigate a data-driven approach to constructing uncertainty sets for robust optimization problems, where the uncertain problem parameters are modeled as random variables whose joint probability distribution is not known. Relying only on independent samples drawn from this distribution, we provide a nonparametric method to estimate uncertainty sets whose  probability mass is guaranteed to approximate a given target mass within a given tolerance with high confidence. The nonparametric   estimators that we consider are  also shown to obey distribution-free finite-sample  performance bounds that imply  their convergence in probability   to the given  target mass. In addition to being efficient to compute, the proposed estimators result in  uncertainty sets that yield computationally tractable robust optimization problems for a large family of constraint functions. 
\end{abstract}


\section{Introduction} \label{sec:introduction}

In this paper, we consider a class of  optimization problems whose feasible regions are defined in terms of \emph{chance constraints} \cite{charnes1959chance}  of the form
\begin{align} \label{eq:chance}
\prob{f(x,  \rand ) \leq 0} \geq \conf.
\end{align}
Here,  $f: \Rset^\ell \times \Rset^d \rightarrow \Rset$ denotes the constraint function,  $x \in \Rset^\ell$ denotes the decision variable, and $\rand$ is an $\Rset^d$-valued random vector that reflects uncertainty in the constraint parameters. The chance constraint \eqref{eq:chance} requires that the  decision variable satisfy the constraint $f(x, \rand) \leq 0$ with probability no smaller than $\conf \in [0,1]$. 

Chance constrained optimization problems are challenging to solve for a variety of reasons. First,  their feasible regions are generally nonconvex \cite{ahmed2008solving}. For example, chance constrained problems have been shown to be NP-hard even in the most basic setting where the constraint function is affine in both the decision variable and uncertain parameters  \cite{luedtke2010integer}. Furthermore, verifying the feasibility of a candidate solution to a chance constrained problem is difficult, because it involves the evaluation of a multivariate integral which can be computationally expensive in high dimensions \cite{ahmed2008solving}. To complicate matters further, the underlying distribution according to which the random vector is distributed may be unknown. 
 
A variety of approaches to the solution of chance constrained problems have been explored in the literature. There are several works that develop exact convex reformulations of chance constrained problems \cite{miller1965chance,prekopa1970probabilistic, prekopa1995stochastic, lagoa1999convexity, lagoa2005probabilistically} and their distributionally robust generalizations \cite{calafiore2006distributionally,yang2016distributionally,zymler2013distributionally}. These reformulations rely on a variety of assumptions, including restrictions on the family of probability distributions according to which the random vector is distributed and the functional form of the chance constraints. Ultimately, many problems of practical interest may fail to satisfy these narrow structural and distributional assumptions.

Because of the rarity of problem instances amenable to exact convex reformulation, there is another line of research focused on the design of approximation methods for chance constrained problems.
One approach that has been explored extensively  involves the construction of explicit convex inner (conservative) approximations to  chance constrained optimization problems \cite{pinter1989deterministic, ben2000robust,nemirovski2006convex,chen2007robust,bertsimas2009constructing, chen2010cvar}. Again, a potential drawback of these methods stems from their reliance (to varying degrees) on knowledge of certain features of the underlying distribution, e.g., support or moment information.
Data-driven approximation methods seek to alleviate the reliance on distributional assumptions that may be overly stringent or difficult to verify in practice. Instead, they utilize data sampled from the underlying distribution. For instance, the sample average approximation method \cite{ahmed2008solving} involves selecting an optimal subset of the sampled data that has empirical probability mass no smaller than the target reliability level $\relLev$. While this approach may be less conservative than other data-driven methods, it gives rise to approximations in the form of mixed-integer optimization problems---which become computationally intractable to solve in the large sample regime.

There is another stream of literature  focusing on data-driven approximations that give rise to tractable convex optimization problems. A particular category of methods uses the given data to construct estimates of the unknown distribution or its lower-order moments \cite{calafiore2006distributionally, delage2010distributionally, ben2013robust, jiang2016data, esfahani2018data}. These estimates, together with  suitably defined confidence regions, give rise to distributionally robust approximations to the original chance constrained program that, in some instances, admit tractable convex reformulations.

Another family of data-driven  methods known as \textit{scenario approximations} have also been widely studied as another tractable alternative to the approximation of chance constrained problems \cite{campi2008exact, calafiore2010random}  and their distributionally robust counterparts \cite{erdougan2006ambiguous, tseng2016random}. 
Specifically, scenario-based methods approximate the  chance constraint \eqref{eq:chance} with $\batNum$ sampled constraints of the form 
\begin{equation} \label{eq:scenario}
f(x, \rand_i) \leq 0, \ i = 1, \dots, n,
\end{equation}
where $\rand_1, \dots, \rand_n$ is an i.i.d. sample drawn from the unknown distribution.
 In addition to being distribution-free, scenario-based methods do not impose restrictions  on the functional form of the constraint function beyond requiring its convexity in the decision variable. However, these methods require a sample size that is at least $O \left(  \frac{\ell}{\vioLev}  \ln\left(\frac{1}{\delta}\right) \right) $ in order to yield a solution that is feasible for the original chance constrained problem with probability at least $1 - \confLev$. While the dependence of this sample size requirement on the problem dimension $\ell$ can be improved for problems that are not fully supported (e.g., \cite{calafiore2010random,schildbach2013randomized,zhang2015sample}), the bound is tight in general.  Therefore, when the number of decision variables is large or when the reliability needed is high, scenario-based approximations can become computationally challenging to solve due to the large number of constraints that must be enforced. 

\subsection{Contribution and Related Work}

 Our approach is predicated on the conservative approximation of  \eqref{eq:chance} in the form of a robust constraint: 
\begin{align} \label{eq:robust}
 f(x,u) \leq 0 \ \ \forall u \in \Ucal,
\end{align}
where the \emph{uncertainty set} $\Ucal \subset \Rset^d$ is constructed from data in  a manner that ensures its satisfaction of the mass requirement $\prob{\rand \in \Ucal} \geq \conf$ with high probability.
In this paper, we provide a nonparametric method to construct such uncertainty sets from data in a manner that guarantees that their probability mass is within a given tolerance of the target mass $\conf$ with high confidence. The uncertainty sets that we propose are trivial to compute, satisfy distribution-free finite-sample statistical guarantees, and give rise to robust  constraints \eqref{eq:robust} that are computationally tractable for a large family of constraint functions.  Additionally, the proposed methodology allows for the explicit representation of nonconvex uncertainty sets as finite unions of convex sets---e.g., as the  union of $m\in \Nset$ balls centered at sampled data points. In this manner, the resulting robust constraint \eqref{eq:robust} can be equivalently reformulated as a finite intersection of $m$ simpler robust constraints, where each robust constraint is defined in terms of a convex uncertainty set. Importantly, the geometry of the underlying convex sets  can be tailored to accommodate the structure of the constraint function to ensure the tractability of the resulting robust constraints.   As another degree of freedom in the parameterization of these sets, the user is free to specify the number of constitutive sets used in the  representation of the uncertainty set, and, hence, the number of robust constraints that must be enforced in the approximation. This provides the user with some degree of control over the size of the resulting robust optimization problem, unlike scenario approximation methods. 

It is important to mention that the probabilistic coverage guarantees accompanying the  class of uncertainty sets proposed in this paper require that the random vector $\uncPar$ have a continuous distribution (cf. Assumption \ref{ass:cont}). The  guarantees accompanying  scenario approximation methods do not require such assumptions.  We also note that, while we suggest a particular class of uncertainty sets in Section \ref{sec:tractability}, the  probabilistic coverage guarantees  provided in Section \ref{sec:guarantees} hold for a more general family of uncertainty sets---e.g., those that satisfy Assumption \ref{ass:cont}.

Our approach is closely related to a class of existing methods that utilize data to construct uncertainty sets that yield tractable robust approximations to chance-constrained optimization problems. For example, Margellos et al. \cite{margellos2014road} develop a data-driven approach to estimate intervals that cover each component of the random vector with high probability, resulting in uncertainty sets that take the form of hyperrectangles. Bertsimas et al. \cite{bertsimas2018data} utilize statistical hypothesis tests to construct convex uncertainty sets from data. However, due to the convex geometry of the uncertainty sets produced by these methods, they may be limited in terms of their ability to accurately describe nonconvex high probability regions that reflect multimodality in the underlying distribution.
Because of the limited expressiveness of these uncertainty sets, they may result in overly conservative approximations to the original chance constraint. Closer to the approach adopted in this paper, Campbell and How \cite{campbell2015bayesian} consider representations in the form of unions of ellipsoids. A potential drawback of their approach, however, stems from the reliance of their theoretical results on the restrictive assumption that the unknown distribution belong to a family of Dirichlet process Gaussian mixtures.
 
The problem of estimating high probability sets from data  also  has similarities to minimum-volume (MV) set estimation problems.
MV-set estimators based on the estimation of density level sets \cite{cadre2006kernel,cuevas1997plug, lei2013distribution,tsybakov1997nonparametric} typically yield nonconvex sets that are not readily expressible as finite unions of convex sets---limiting their tractability from a robust optimization perspective. More broadly, MV-set estimation methods include  those based on nonparametric set estimation \cite{devroye1980detection,ba2000set}, empirical quantile-based estimation \cite{polonik1997minimum, scott2006learning},  classification via support vector machines  \cite{scholkopf2001estimating, muandet2012learning}, and network flow-based methods \cite{huo2004network}. With the exception of  \cite{scott2006learning}, the methods presented in these papers lack explicit finite-sample guarantees on their performance. A practical limitation of the algorithms proposed by Scott et al. \cite{scott2006learning}, however, is that the computation of their MV-set estimates is intractable in high dimensions.

 \subsection{Organization} 
The remainder of the paper is organized as follows. In Section \ref{sec:formulation}, we introduce a data-driven method to construct uncertainty sets with a given probability mass. In Section \ref{sec:guarantees}, we provide finite-sample statistical guarantees on the performance of the proposed class of estimators. Section \ref{sec:experiments} illustrates the behavior of the estimators with several numerical experiments, and Section \ref{sec:conclusion} concludes the paper.

\subsection{Notation} 
 We employ the following notational conventions throughout the paper.
 Let $\Rset$ and $\Zset$ denote the sets of real numbers and  integers, respectively. Given a positive integer $n \in \Zset$, we let $[n] := \{1, \dots, n\}$ denote the set of the first $n$ positive integers.  Given a real number $x \in \Rset$, we denote its \emph{ceiling} by $\lceil x \rceil : = \min\{ n \in \Zset \ | \ n \geq x \}$. 
Throughout, we  use boldface symbols to denote random variables, and non-boldface symbols to denote particular values in the range of a random variable and other deterministic quantities. We  use $\prob{A}$ to denote the probability of an event $A$.

\section{Uncertainty Set Estimation} \label{sec:formulation} 
In this section, we provide a data-driven method to construct uncertainty sets  whose probability mass is both close to and no smaller than a given target mass $\conf \in (0,1)$ with high probability, based only on a finite random sample drawn from the unknown distribution of interest.
In addition to the mass requirement, we seek representations that result in computationally tractable robust constraints. 

The remainder of this section is organized as follows. In Sec. \ref{sec:level}, we provide a general  approach to learning  $\conf$-probability  sets via  the estimation of sublevel sets for a general class of continuous functions. In Sec. \ref{sec:tractability}, we introduce a particular class of functions whose level sets (expressible as finite unions of convex sets) yield tractable robust constraints for several important families of constraint functions. We close with a discussion on the connections between the proposed family of approximations and scenario-based approximation methods.

\subsection{Estimation via Level Sets} \label{sec:level}
With the aforementioned objectives in mind,  we consider a family of uncertainty sets  defined as  $r$-sublevel sets\footnote{The decision to specify \eqref{eq:estimator} in terms of sublevel sets of $\phi$, as opposed to superlevel sets, is made for notational convenience.} of a  given continuous function $\phi: \Rset^d \rightarrow \Rset$, i.e.,
\begin{align} \label{eq:estimator}
\lev{r} := \{ u \in \Rset^d \ | \ \phi(u) \leq r\}.
\end{align}

The proposed family of uncertainty sets \eqref{eq:estimator} is parameterized  by the  user-specified inputs: the \emph{shape function} $\phi$ and \emph{level} $r$. Intuitively, the shape function $\phi$ should be chosen so that it is smallest over those regions with the greatest concentration of probability mass, in order to limit the volume (Lebesgue measure) of the resulting sublevel set. While this might suggest an approach to  specifying  $\phi$ in terms of a  density estimate based on an i.i.d. sample drawn from the unknown distribution, the particular functional form  of $\phi$ must also yield uncertainty (level) sets that are computationally tractable from a robust optimization perspective---a condition that most nonparametric (kernel) density estimators fail to satisfy. 
In Sec. \ref{sec:tractability}, we suggest a particular functional form for the shape function $\phi$ that is expressive enough to represent compact sets of arbitrary form, while ensuring the tractability of the robust constraints that it yields for a large family of constraint functions.  It is important to note that, in the interest of lightening notation, we have omitted the potential dependence that the function $\phi$ may have on data drawn from the unknown distribution of interest.

For the remainder of this section, we treat the shape function $\phi$ as being fixed (implicitly conditioning all probabilistic statements on the data on which it is based), and focus our discussion on the role of the level $r$ in controlling the behavior  of the sublevel sets of $\phi$. Clearly, the level $r \in \Rset$ directly influences the volume and probability mass of the  sublevel sets of $\phi$, with larger levels  resulting in  sublevel sets of greater volume. Thus, with the aim of  limiting the conservatism of the robust approximations \eqref{eq:robust} induced by the proposed family of uncertainty sets, we are interested in characterizing the smallest level $r$ such that the target mass requirement is satisfied. We refer to this as the $\conf$-\emph{covering level}, which we formally define as follows.  
 
 \begin{defn}[$\conf$-covering level] Let $\conf \in (0,1)$. The $\conf$-\emph{covering level} is defined as
\begin{equation*}
\rho(\conf) := \inf \left\{r  \in \Rset \ | \ \prob{\rand \in \lev{r}} \geq \conf \right\}.
\end{equation*}
\end{defn}
\vspace{.05in}

As one of the objectives of this paper, we are interested in constructing an estimator for the  $\conf$-covering level $\rho(\conf)$ based only on an  i.i.d. training sample $\rand_1, \dots, \rand_n$ drawn from the unknown distribution. Before introducing the particular family of estimators that we consider in this paper, it will be informative to express the $\conf$-covering level in terms of the  quantile 
$$\rho(\conf) = F^{-1}(\alpha),$$
where $F^{-1}(\alpha) := \inf\{z \in \Rset \ | \ \prob{\bzeta \leq z} \geq \conf \}$ denotes the  $\conf$-quantile of the transformed random variable $\bzeta := \phi(\rand)$. This  reformulation suggests a natural estimator for  $\rho(\conf)$ in the form of an empirical quantile
\begin{equation} \label{eq:quantile}
F_n^{-1}(\conf_n) := \inf\{  z \in \Rset \ | \  F_n(z) \geq \conf_n \}.
\end{equation}
Here, $\conf_n \in (0,1)$ is a sequence of probability levels, and   $F_n(z) := (1/n) \sum_{i=1}^n \mathbf{1} \{ \bzeta_i \leq z\}$ denotes the empirical distribution function associated with the transformed training sample $$\bzeta_i := \phi(\rand_i)$$ for $i = 1, \dots, n$.

In Section \ref{sec:guarantees}, we investigate the role of $\conf_n$ in controlling the limiting behavior and rate of convergence of the proposed empirical quantile to the $\conf$-covering level. In particular, we establish minimum training sample size requirements on $n$ ensuring that the probability mass covered by the resulting uncertainty set  $\lev{F_n^{-1}(\conf_n)}$ is within a given tolerance of (and no smaller  than) the given target mass $\conf$ with  high probability.

\subsection{Shape Functions $\phi$ Compatible with Robust Optimization}\label{sec:tractability}

 The functional form of the shape function $\phi$ (and the geometry of its sublevel sets) will play a critical role in governing the tractability of the family of robust constraints that it gives rise to.  With this in mind, we suggest a class of shape functions  whose sublevel sets are rich enough to represent compact sets of arbitrary form, while being simple enough to ensure the tractability of the resulting robust constraints. Specifically, we propose a particular class of shape functions whose sublevel sets admit explicit representations as finite unions of regular convex sets, e.g., $p$-norm balls. In this manner, the resulting robust constraint \eqref{eq:robust} can be equivalently reformulated as a finite intersection of simpler robust constraints, where each robust constraint is defined in terms of a convex uncertainty set. Being expressed as norm balls, the particular geometry of the underlying convex sets (e.g., $\ell_1$ versus $\ell_2$) can be tailored to accommodate the structure of the given constraint function to facilitate the tractability of the corresponding robust constraint.
In addition to determining the number of constraints, the number of sets in the union will also influence the volume of the uncertainty sets that we learn, and hence the conservatism of the resulting robust constraint approximations. 
Therefore, in treating the number of sets in the approximation as a user-specified parameter, one can tradeoff the computational complexity of the resulting optimization problem against the quality (conservatism) of its solutions.

With these design criteria in mind, we now introduce the specific class of estimators we consider.
Given an i.i.d. sample $\bar{\uncPar}_{1}, \dots, \bar{\uncPar}_\batSiz$ drawn from the unknown distribution of interest,  define
\begin{align} \label{eq:ballest}
	\phi \left(  u  \right) := \min_{i \in \left[ \batSiz  \right]} \hnorm{ u - \bar{\uncPar}_{i}  },
\end{align}
where $\hnorm{\cdot}$ denotes the $p$-norm ($p \geq 1$) on $\mathbb{R}^d$.  The specification of the shape function according to \eqref{eq:ballest} gives rise to a class of uncertainty sets that are defined as a union of balls with a common radius, where each ball is centered at a randomly sampled point
\begin{equation}\label{eqn:unc_set_def}
	\uncSet_\phi \left( r \right) = \bigcup_{\batSizMem=1}^\batSiz B \left( \bar{\uncPar}_{\batSizMem}, r   \right).
\end{equation} 
Here $B \left( \uncRel, r   \right) := \left\{ v  \in \mathbb{R}^d \;|\;   \hnorm{u-v} \leq r   \right\}$ denotes the closed ball of radius $r \geq 0$ centered at $\uncRel \in \R^d$. The parameter $r$ plays an analogous role to that of bandwidth parameters in the context of nonparametric kernel density estimation.\footnote{Nonparametric estimators defined as the union of balls centered on the sampled data have been previously studied in the context of distribution support estimation \cite{ba2000set, devroye1980detection}.}

It is important to note that we have augmented the notation used to denote the sample $(\bar{\uncPar}_{1}, \dots, \bar{\uncPar}_\batSiz)$ (which we refer to as the \emph{shape sample}), in order to distinguish its specification from the  \emph{training sample}  $(\rand_1, \dots, \rand_n)$ that is subsequently used to estimate the $\conf$-covering level. The two samples, which are assumed to be independent from one another, play distinct roles in the construction of the uncertainty set. The shape sample determines the positions of the balls comprising the uncertainty set, while the training sample is used to calibrate the radius $r$ of each ball in a manner guaranteeing that the probability mass covered by the resulting uncertainty set is within a given tolerance of the target probability mass with high confidence.

Given an uncertainty set $\lev{r}$ defined according to \eqref{eqn:unc_set_def}, one can equivalently reformulate the robust constraint \eqref{eq:robust} as the intersection of $\batSiz$  robust constraints of the form
\begin{align} \label{eq:robust_reform}
f \left( x, \uncRel   \right) \leq 0 \quad \forall \uncRel \in B \left( \bar{\uncPar}_{\batSizMem}, r   \right), \quad \batSizMem = 1, \dots, \batSiz.
\end{align}
Robust constraints defined in terms of uncertainty sets of this form result in computationally tractable optimization problems for a large family of constraint functions. For example, if the given  constraint function is a \emph{bi-affine function} of the form $f(x,u) = x^\top u - b$,
then it is straightforward to show that the collection of robust constraints  specified in  \eqref{eq:robust_reform} can be equivalently reformulated as convex constraints given by
\begin{align}  \label{eq:dual}
x^\top \bar{\uncPar}_i  + r \hnorm{x}_* \leq b, \quad \batSizMem = 1, \dots, \batSiz,
\end{align}
where $\hnorm{x}_* := \sup \{ x^\top y  \ | \  \hnorm{y} \leq 1\}$ denotes the \emph{dual norm} associated with a given $p$-norm  $\hnorm{\cdot}$.

We refer the reader to \cite{bertsimas2006tractable} and \cite{ben2009robust} for a more comprehensive discussion surrounding the different families of robust constraints that admit  equivalent reformulations as tractable convex constraints. 

\subsection{Comparison to Scenario Approximation Methods}

The approach considered here gives rise to optimization problems that are similar in structure to those based on scenario approximations of chance constrained problems and their distributionally robust generalizations. Specifically, when the radius $r$ of the balls composing the uncertainty sets is equal to zero, we recover the standard scenario approximations to chance constrained problems such as those presented in \cite{campi2008exact, calafiore2010random}. When $r > 0$, we recover the robust scenario approximation to distributionally robust chance constrained problems studied in \cite{erdougan2006ambiguous, tseng2016random}. Interestingly, the latter connection suggests that the class of approximations considered in this paper also possesses an intrinsic distributional robustness---a point that merits further examination as part of future research.

Although the class of approximations that we propose is similar in structure to those based on scenario approximation methods, our method offers an important computational advantage in terms of the number of constraints that must be enforced in the resulting approximation to the original chance constraint. 
Specifically, while scenario approximation methods result in approximations based on a number of sampled constraints that is at least  $O(\frac{\ell}{1 - \alpha} \ln \left(\frac{1}{\delta} \right))$, the method proposed in this paper treats the number of constraints $m$ in the resulting robust approximation as a parameter that can be set by the user. Importantly, our theoretical guarantees hold for any value of $m$, allowing the user to determine the number of constraints that must be enforced in the resulting approximation to the original chance constraint.  

It is also important to note that the class of approximations proposed in this paper is accompanied by  certain limitations,  as compared to  scenario approximation methods. In particular,  while scenario approximations \eqref{eq:scenario} only require that the constraint function be convex in the decision variable,  the tractability of our approximations  \eqref{eq:robust_reform} requires  additional restrictions on the class of constraint functions that can be efficiently optimized over---e.g., that the constraint functions  are also concave with respect to the uncertain variable.

 \section{Finite-Sample Statistical Guarantees}\label{sec:guarantees}
 We now establish a bound on the rate at which the probability mass of the proposed class of uncertainty sets converges to the target mass as a function of the training sample size $n$. It is important to note that the theoretical results presented in this section hold 
 not only for the class of uncertainty sets introduced in Section \ref{sec:tractability}, but more generally for any family of uncertainty sets  induced by shape functions that satisfy Assumption \ref{ass:cont}. To lighten the notation in the following discussion, we denote the probability mass of the uncertainty set $\lev{r}$ conditioned on a level $r$ by
\begin{align} \label{eq:probmass}
 \pi_\phi (r) := \prob{\rand \in \lev{r}}.
\end{align} 
Throughout this section, we treat the shape function $\phi$ as being fixed and given. Accordingly,  all  subsequent probabilities that we  state in this section should be  interpreted as being conditioned on the given shape function $\phi$. 
Before stating the main result of this section, we introduce a technical assumption that is assumed to hold throughout the remainder of the paper.

 \begin{ass} \label{ass:cont} Given a fixed shape function $\phi: \Rset^d \rightarrow \Rset$, the transformed random variable $\bzeta = \phi(\rand)$  is assumed  to have a continuous distribution.
\end{ass}
\vspace{.05in}
In particular, Assumption \ref{ass:cont} is satisfied by any shape function  $\phi$ constructed  according to \eqref{eq:ballest} if the underlying random vector $\rand$ is continuous. With Assumption \ref{ass:cont} in hand, we have the following result.

\begin{thm} \label{thm:bound} Let $\varepsilon \in (0, 1- \conf)$ and  $\conf_n \in (\conf, \conf + \varepsilon)$ for all $n \geq 1$. It follows that \vspace{-.05in}
\begin{equation} \label{eq:error1}
\prob{\pi_\phi (F_n^{-1}( \conf_n)) < \conf } \leq \exp\left( \frac{-n(\conf_n - \conf)^2}{2(1- \conf)} \right)
\end{equation}
and \vspace{-.05in} 
\begin{align} \label{eq:error2}
\nonumber \prob{ \pi_\phi (F_n^{-1}( \conf_n)) > \conf + \varepsilon }  \qquad \qquad \qquad \qquad  \\
\qquad \qquad \qquad \qquad \leq  \exp\left( \frac{-n(\conf + \varepsilon - \conf_n)^2}{2(\conf + \varepsilon)}    \right).
\end{align}
\end{thm}
\vspace{.08in}

It follows that, for any fixed tolerance $\varepsilon \in (0, 1-\conf)$, the probability that the mass covered by the uncertainty set $\lev{F_n^{-1}(\conf_n)}$  is within $\varepsilon$ of (and no less than) the target mass $\conf$ approaches one at an exponential rate in the sample size. One may also use this result to establish conditions on the sequence $\conf_n$ guaranteeing that $\pi_\phi (F_n^{-1}(\conf_n)) \rightarrow \conf$  in probability as the sample size $n$ tends to infinity.

We note that the proof of Theorem \ref{thm:bound} utilizes a number of standard arguments that have been used to establish several related results in the literature \cite{khargonekar1996randomized,chen1998order,alamo2010sample}. The crux of the main argument involves the reformulation of the error probabilities \eqref{eq:error1} and \eqref{eq:error2} as binomial tail probabilities, which are then bounded from above using a variant of the Chernoff bound. 

\begin{proof} All probabilities stated in this proof are assumed to be conditioned on the given shape function $\phi$. We first  prove   inequality \eqref{eq:error1}.
First notice that, conditioned on a level $r$, it holds that $$\pi_\phi(r) = \prob{\phi(\rand) \leq r} = F(r),$$ where $F$ denotes the cumulative distribution function of the transformed random variable $\bzeta := \phi(\rand)$. This identity, combined with the assumed continuity of the distribution $F$, implies that
\begin{align*}
\prob{\pi_\phi(F_n^{-1}( \conf_n)) < \conf}  = \prob{F_n^{-1}( \conf_n) < F^{-1}(\conf)  }.
\end{align*}
 The empirical quantile function can be expressed in terms of the \emph{order statistics} of the transformed training sample, which we denote by $\bzeta_{(1)} \leq \bzeta_{(2)} \leq  \cdots \leq \bzeta_{(n)}$. Specifically,  it holds that 
 \begin{align} \label{eq:score}
  F_n^{-1}(\gamma) = \bzeta_{(\lceil n \gamma \rceil)}
 \end{align}
for all  $\gamma \in (0,1)$.  Importantly, the cumulative distribution function of the $k$-th order statistic can be expressed according to the upper tail of a binomial:
\begin{align} \label{eq:orderbin}
\prob{\bzeta_{(k)} \leq z   } = \sum_{i =k}^n {n\choose i}F(z)^i(1- F(z))^{n-i}
\end{align}
for each $k = 1, \dots,n$ \cite{david2004order}. It follows that
  \begin{align}
 \nonumber & \prob{F_n^{-1}( \conf_n) < F^{-1}(\conf)  } \\
 \nonumber & \hspace{1in}= \prob{ \bzeta_{(\lceil n\conf_n \rceil)} < F^{-1}(\conf)}\\
\label{eq:upper} &  \hspace{1in}= \sum_{i =  \lceil n\conf_n \rceil }^n {n\choose i} \conf^i(1- \conf)^{n-i}\\
  \label{eq:binomial} &  \hspace{1in}= \sum_{i =  0}^{n- \lceil n\conf_n \rceil } {n\choose i} (1-\conf)^i \conf^{n-i}.
 \end{align}
 The first equality follows from \eqref{eq:score}. The second equality follows from \eqref{eq:orderbin} and  Assumption \ref{ass:cont}, which implies that $F(F^{-1}(\conf)) = \conf$. The third equality stems from an equivalent reformulation of the  upper binomial tail \eqref{eq:upper} as a lower binomial tail. One can  bound \eqref{eq:binomial} from above using the following well-known upper bound on the lower tail of a binomial \cite{upfal2005probability}[Theorem 4.5], which is  a direct consequence of the Chernoff bound.
 
 \begin{thm} \label{thm:cher} Let $\boldsymbol{\xi}$ be a binomial random variable with parameters $n \in \Nset$ and $p \in [0,1]$. For $k \leq np$, it holds that $$  \prob{\boldsymbol{\xi} \leq k} \leq \exp\left( \frac{-(np - k)^2}{2np} \right).  $$
 \end{thm}
 \vspace{.08in}
Inequality  \eqref{eq:error1} follows from an application of Theorem \ref{thm:cher} to \eqref{eq:binomial}, where we use the fact that $n- \lceil n\conf_n \rceil \leq n(1 - \conf_n)$.
 
 The proof of inequality \eqref{eq:error2} is analogous in nature to the proof of \eqref{eq:error1}. Using similar arguments, it is possible to show that 
 \begin{align}
 \nonumber    & \prob{ \pi_\phi (F_n^{-1}( \conf_n)) > \conf + \varepsilon  }   \\
    \label{eq:binomial2} & \hspace{.5in}  = \sum_{i =  0}^{\lceil n\conf_n \rceil -1 } {n\choose i} (\conf + \varepsilon)^i (1- \conf - \varepsilon)^{n-i}.
 \end{align}
Inequality  \eqref{eq:error2} follows from an application of Theorem \ref{thm:cher} to \eqref{eq:binomial2}, where we use the fact that $\lceil n\conf_n \rceil - 1\leq n\conf_n$. \end{proof}

It is also possible to use Theorem \ref{thm:bound} to characterize a distribution-free  bound on the sample size requirement ensuring that the  probability mass covered by the uncertainty set satisfies the given tolerance with a given confidence $1- \delta$. We state the following corollary without proof, as it is an immediate consequence of Theorem \ref{thm:bound}.
\begin{cor} Let $\delta \in (0,1)$,  $\varepsilon \in (0, 1- \conf)$,  and $\lambda \in (0,1)$.  Set $\conf_n = \conf + \lambda \varepsilon$ for all $n \geq 1$. If 
\begin{equation} \label{eq:samplesize}
n \geq c(\lambda, \conf, \varepsilon)   \left(\frac{2}{\varepsilon^2}\right)\ln \left(\frac{2}{\delta}\right),
\end{equation}
where $c(\lambda, \conf, \varepsilon) := \max\{(1-\conf)/\lambda^2, (\conf + \varepsilon)/(1-\lambda)^2 \},$ then
 $$\prob{ \conf \leq \pi_\phi (F_n^{-1}( \conf_n)) \leq \conf + \varepsilon   }  \geq 1- \delta.$$
\end{cor}
\vspace{.08in}

It is important to note that the sample size requirement \eqref{eq:samplesize} is dimension-free in that it does not depend on the dimension of the random vector or that of the decision variable. Additionally, for $\alpha \in (1/2,1)$, it is straightforward to show that the value of $\lambda$ which minimizes the sample size requirement \eqref{eq:samplesize} is given by
\begin{equation}\label{eq:optimalLambda}
    \lambda^{*} = \frac{1-\relLev - \sqrt{\left(1-\relLev \right) \left(\relLev + \buf \right)}}{1-2 \relLev - \buf}.
\end{equation}

 We also remark that it is possible to improve upon the sample size requirement \eqref{eq:samplesize} through a refinement of the upper bounds on the lower binomial tails  \eqref{eq:binomial}-\eqref{eq:binomial2} using arguments analogous to those in \cite{alamo2010sample}. In particular, the dependence on the tolerance parameter $\buf$ can be improved to $O(1/\varepsilon)$.

\begin{figure*}
	
		\begin{subfigure}[b]{.19\linewidth}
		\centering \includegraphics[width=1\linewidth]{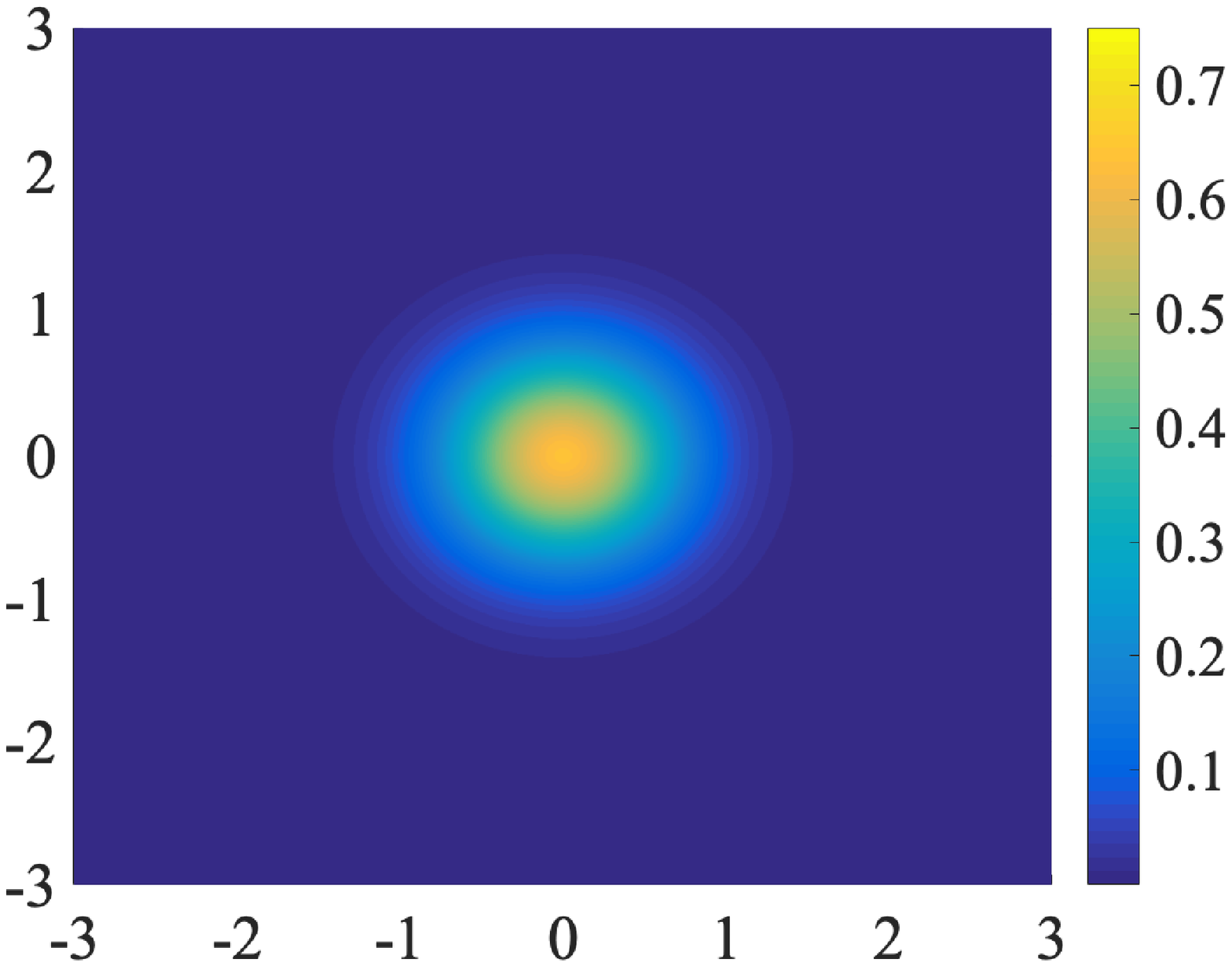}
		\caption{True density}
		\label{fig:density3}
	\end{subfigure}
	\begin{subfigure}[b]{.19\linewidth}
		\centering \includegraphics[width=1\linewidth]{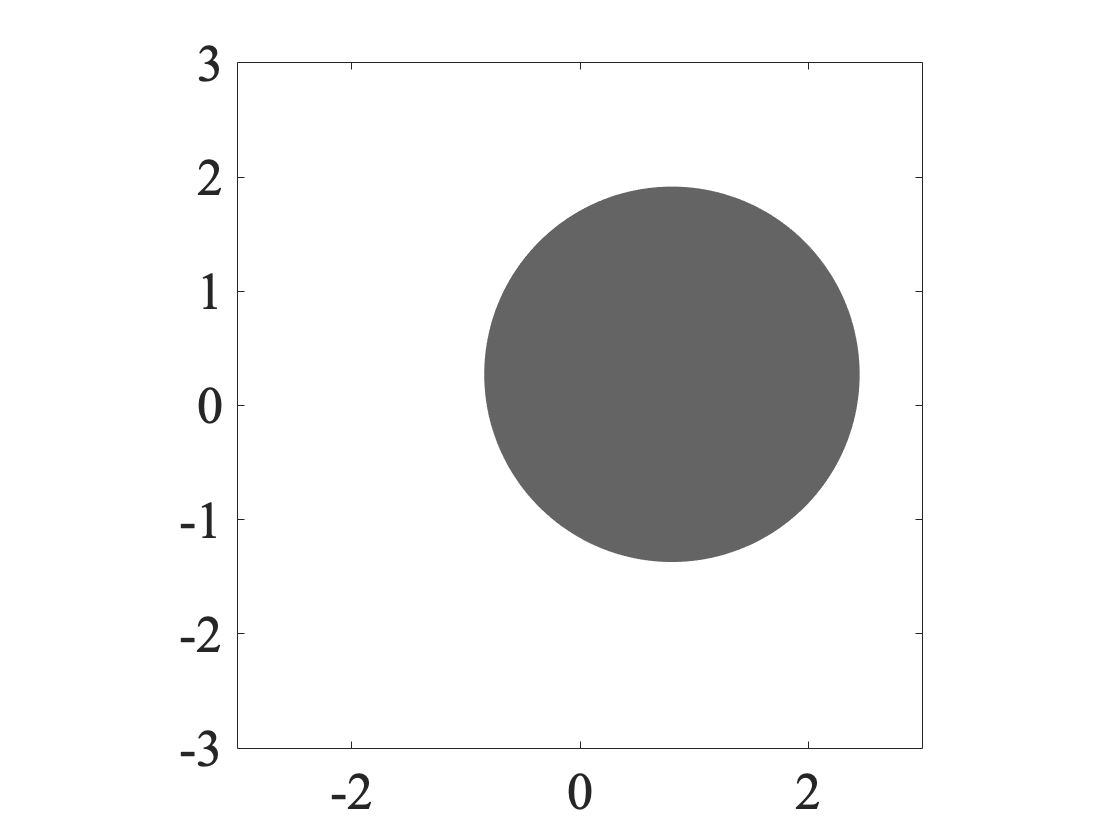}
		\caption{$\batSiz = 1$}
		\label{fig:batSizShapes13}
	\end{subfigure}
	\begin{subfigure}[b]{.19\linewidth}
		\centering \includegraphics[width=1\linewidth]{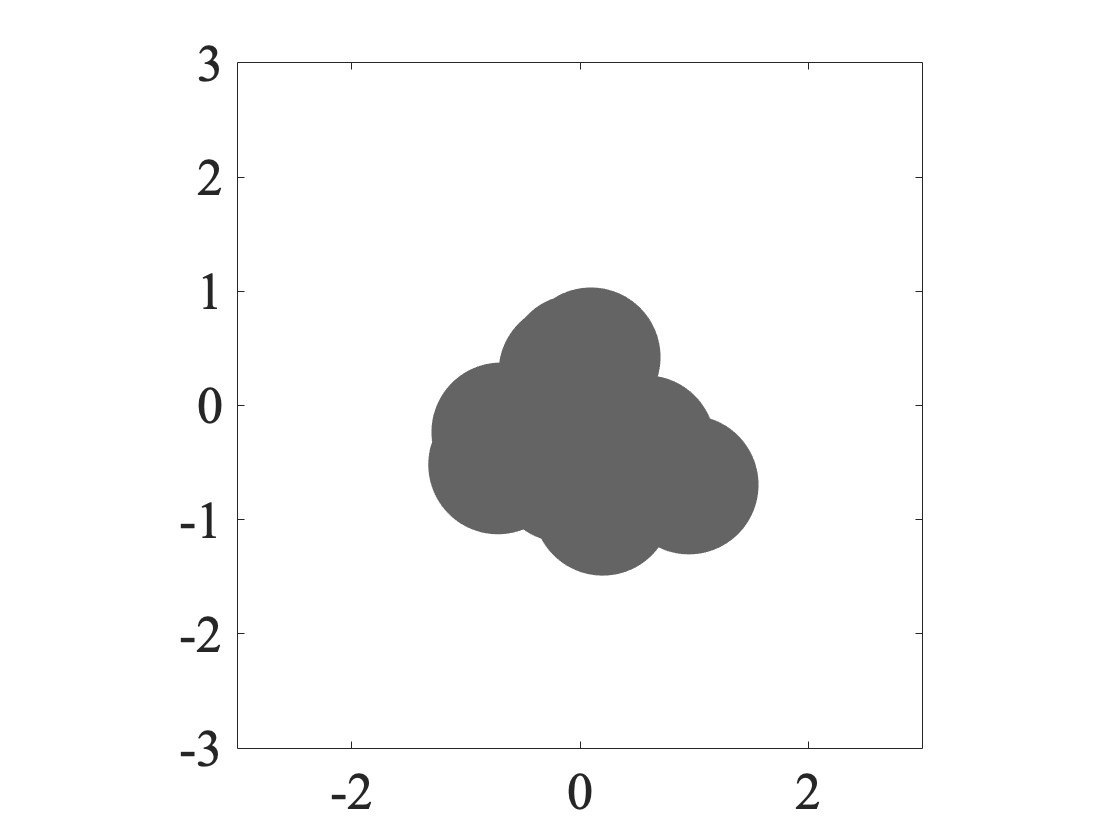}
		\caption{$\batSiz = 10$}
		\label{fig:batSizShapes103}
	\end{subfigure}
	\begin{subfigure}[b]{.19\linewidth}
		\centering \includegraphics[width=1\linewidth]{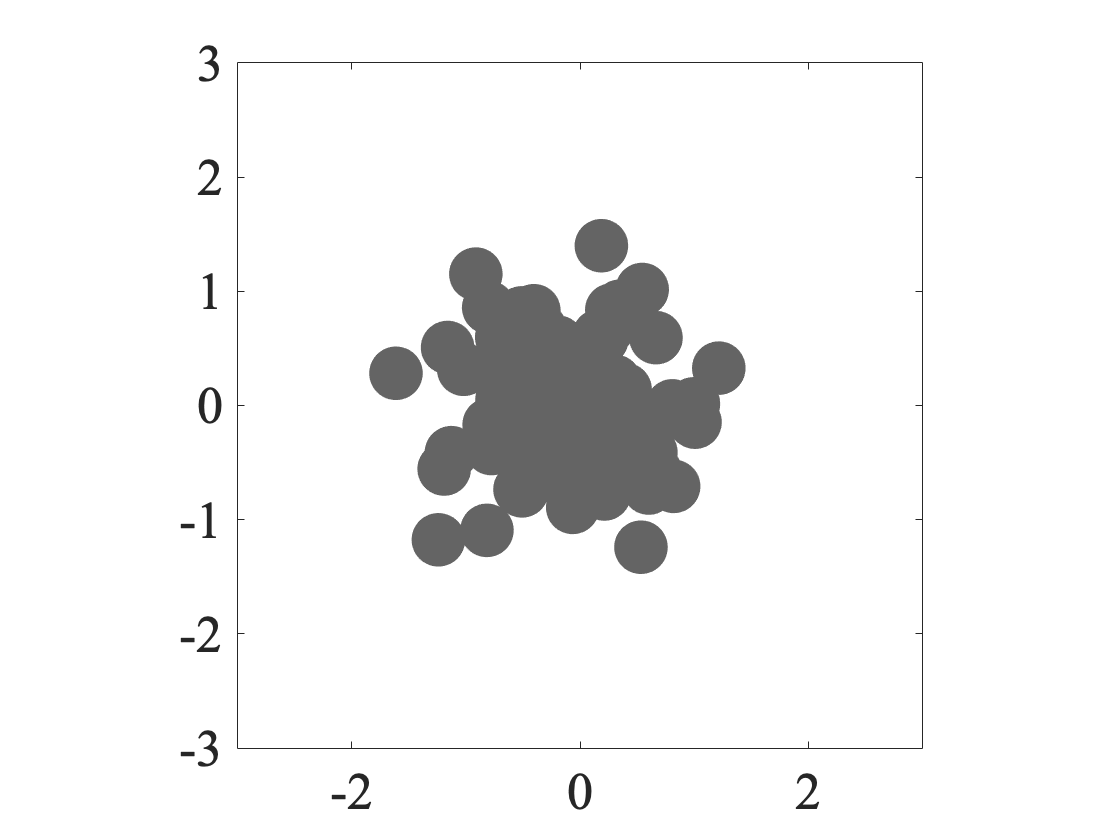}
		\caption{$\batSiz = 100$}
		\label{fig:batSizShapes1003}
	\end{subfigure}
	\begin{subfigure}[b]{.19\linewidth}
		\centering \includegraphics[width=1\linewidth]{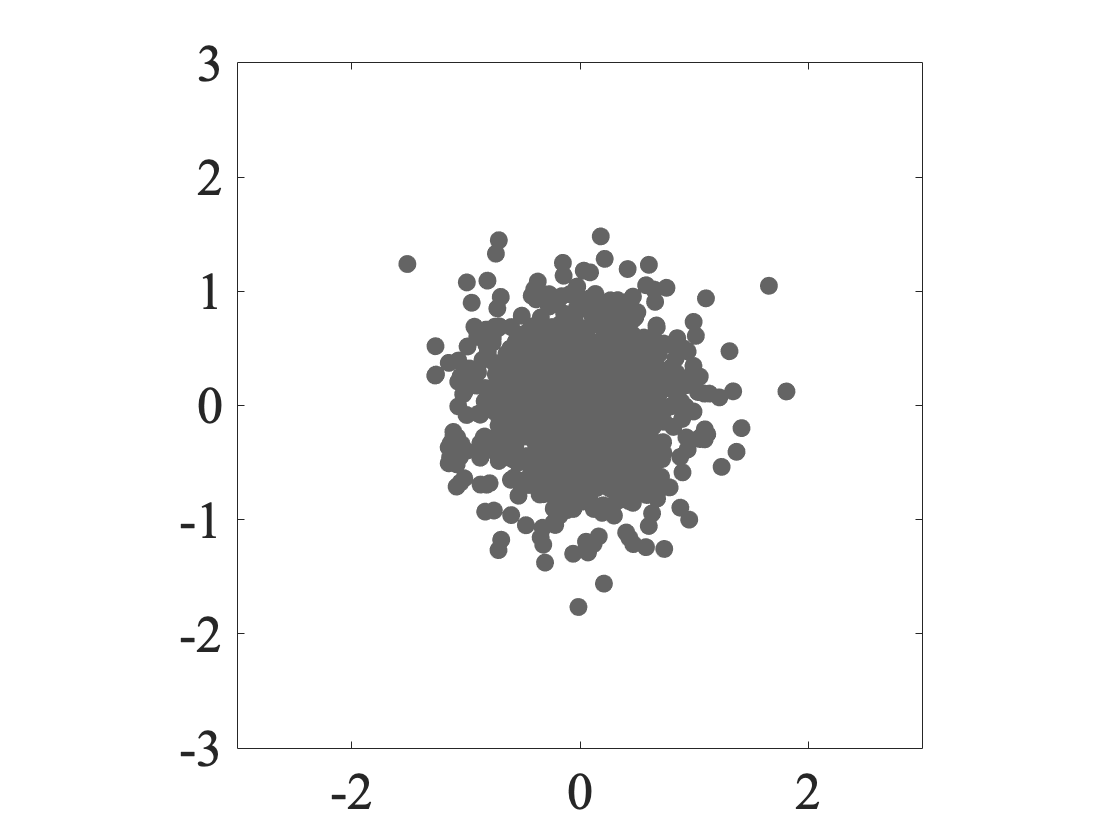}
		\caption{$\batSiz = 1000$}
		\label{fig:batSizShapes10003}
	\end{subfigure}
	
	\begin{subfigure}[b]{.19\linewidth}
		\centering \includegraphics[width=1\linewidth]{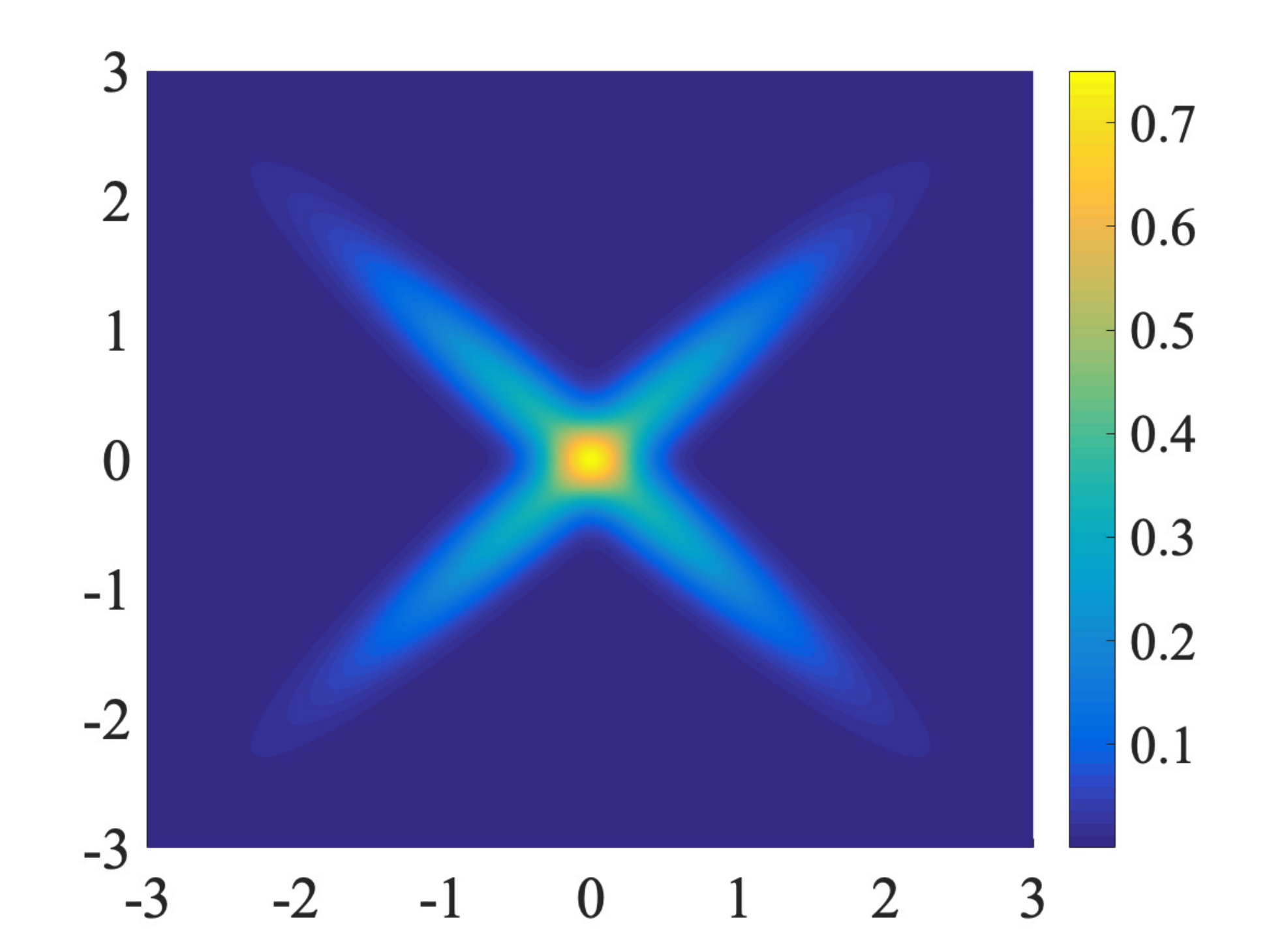}
				\caption{True density}
		\label{fig:density}
	\end{subfigure}
\begin{subfigure}[b]{.19\linewidth}
		\centering \includegraphics[width=1\linewidth]{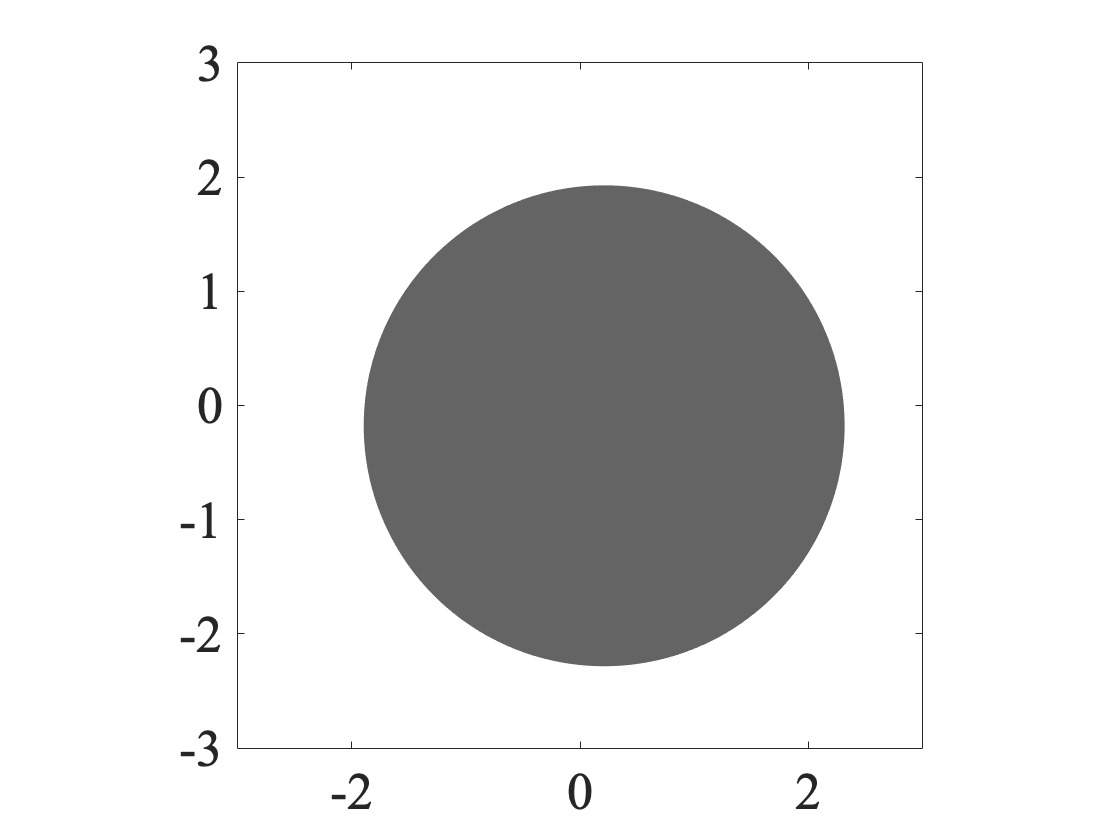}
		\caption{$\batSiz = 1$}
		\label{fig:batSizShapes1}
\end{subfigure}
\begin{subfigure}[b]{.19\linewidth}
	\centering \includegraphics[width=1\linewidth]{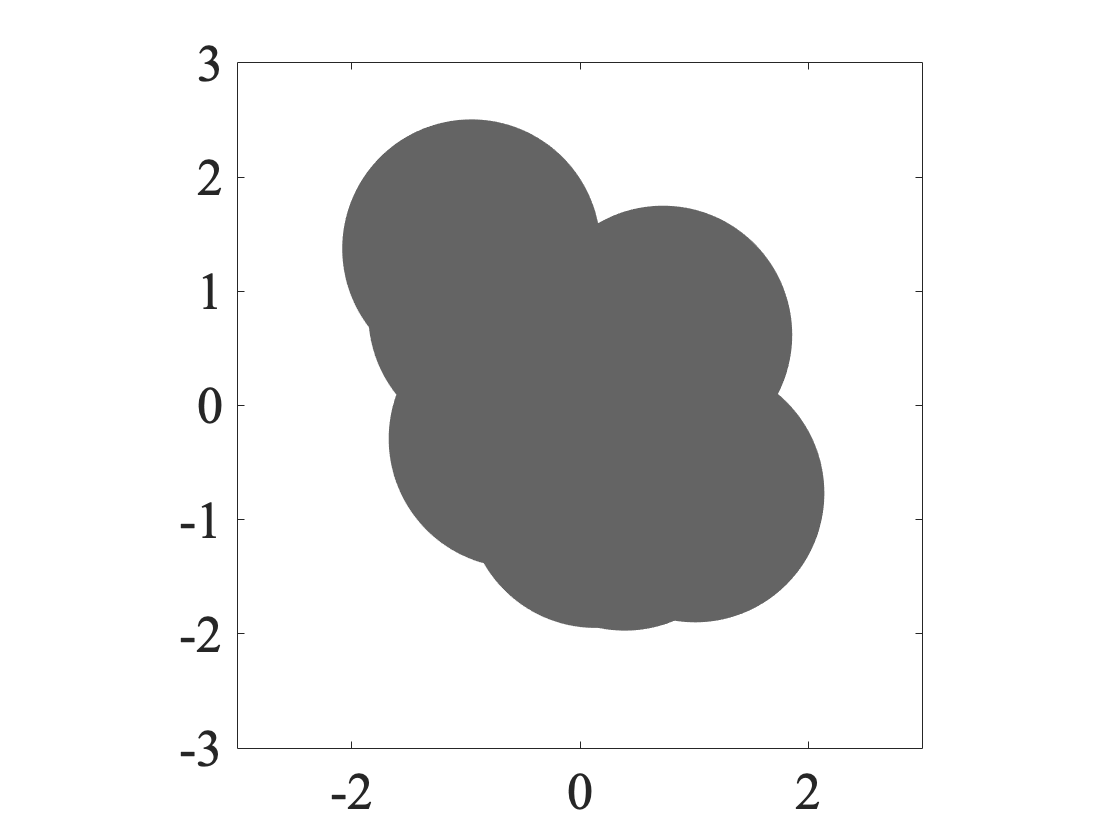}
		\caption{$\batSiz = 10$}
		\label{fig:batSizShapes10}
\end{subfigure}
\begin{subfigure}[b]{.19\linewidth}
	\centering \includegraphics[width=1\linewidth]{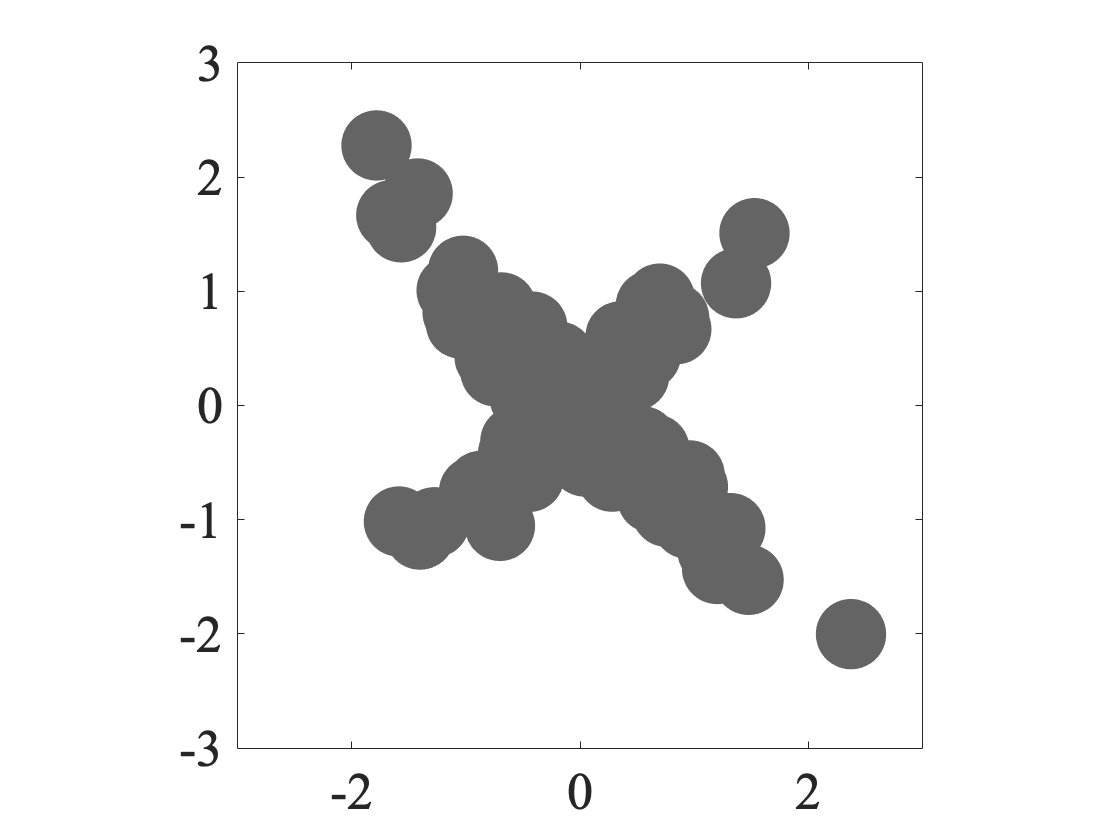}
		\caption{$\batSiz = 100$}
			\label{fig:batSizShapes100}
\end{subfigure}
\begin{subfigure}[b]{.19\linewidth}
	\centering \includegraphics[width=1\linewidth]{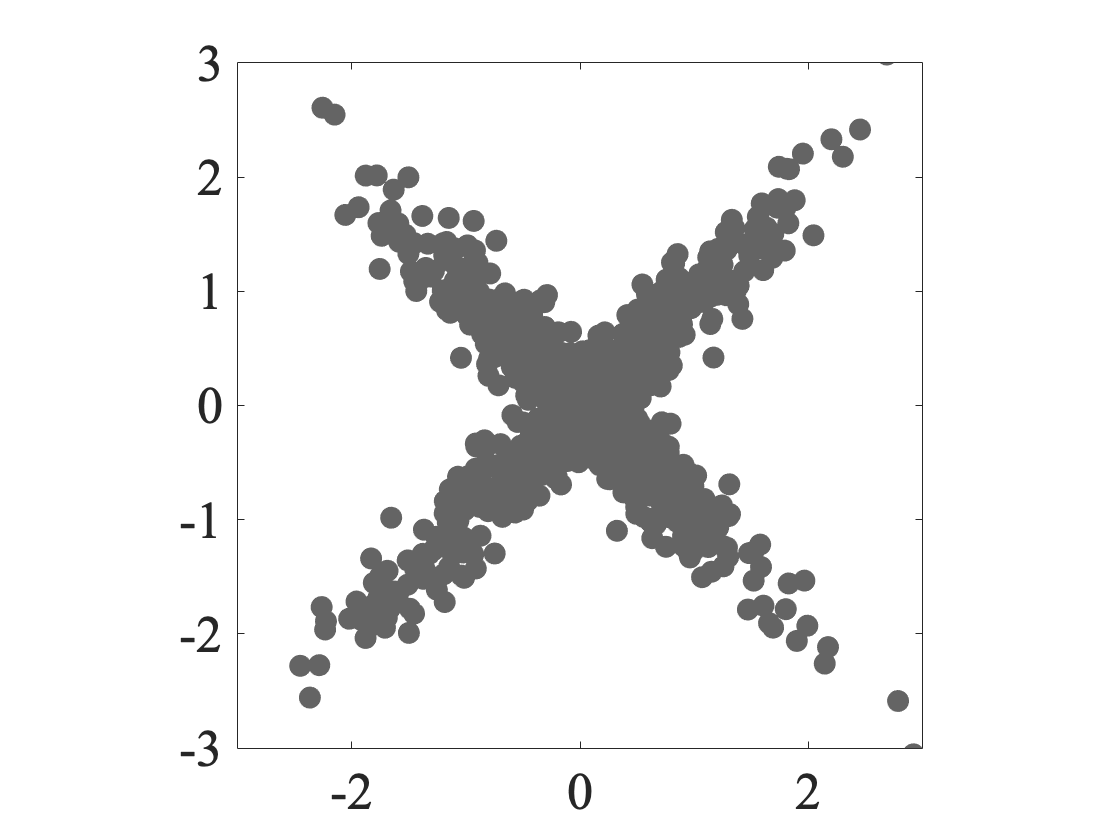}
		\caption{$\batSiz = 1000$}
		\label{fig:batSizShapes1000}
\end{subfigure}

	\begin{subfigure}[b]{.19\linewidth}
		\centering \includegraphics[width=1\linewidth]{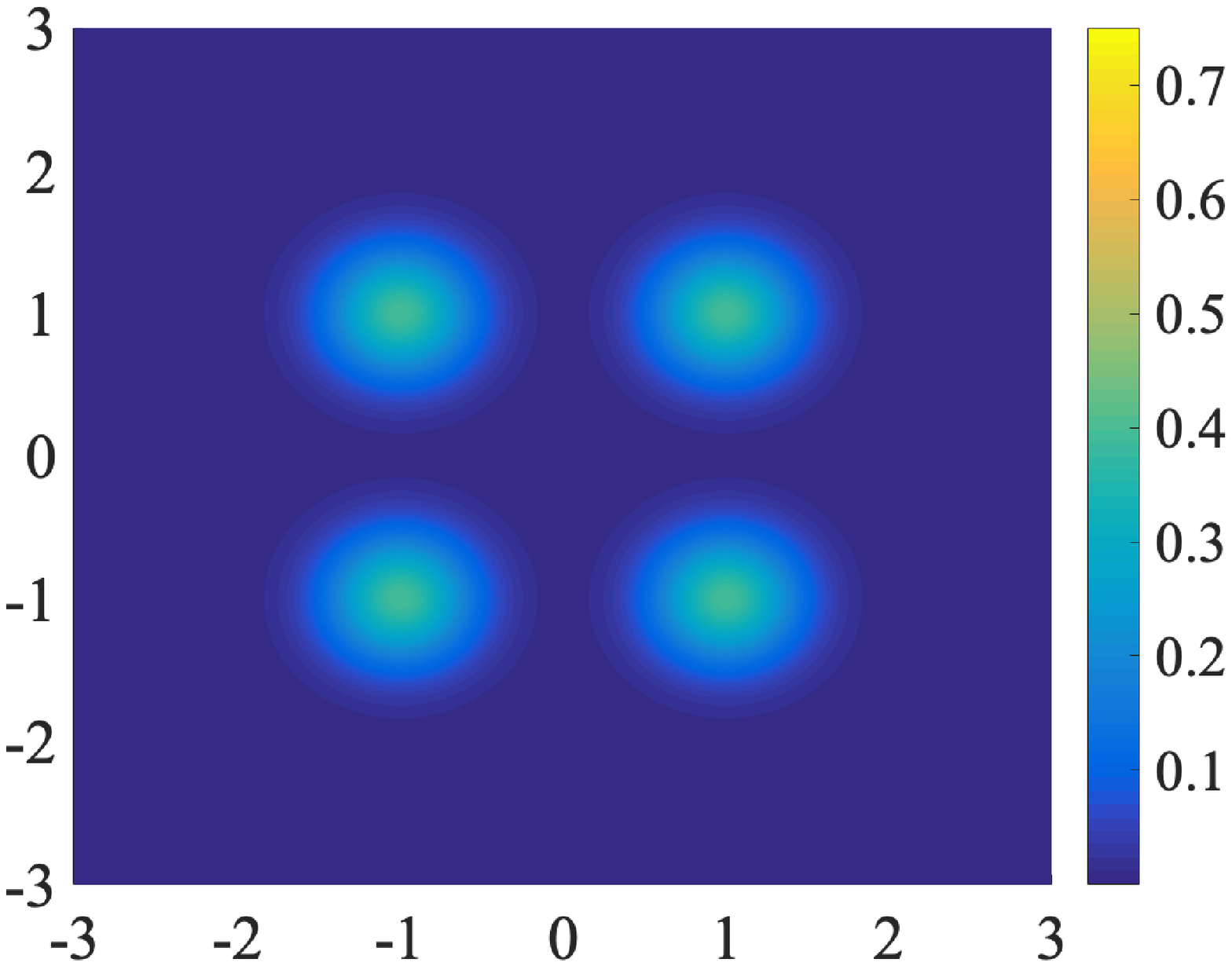}
		\caption{True density}
		\label{fig:density2}
	\end{subfigure}
	\begin{subfigure}[b]{.19\linewidth}
		\centering \includegraphics[width=1\linewidth]{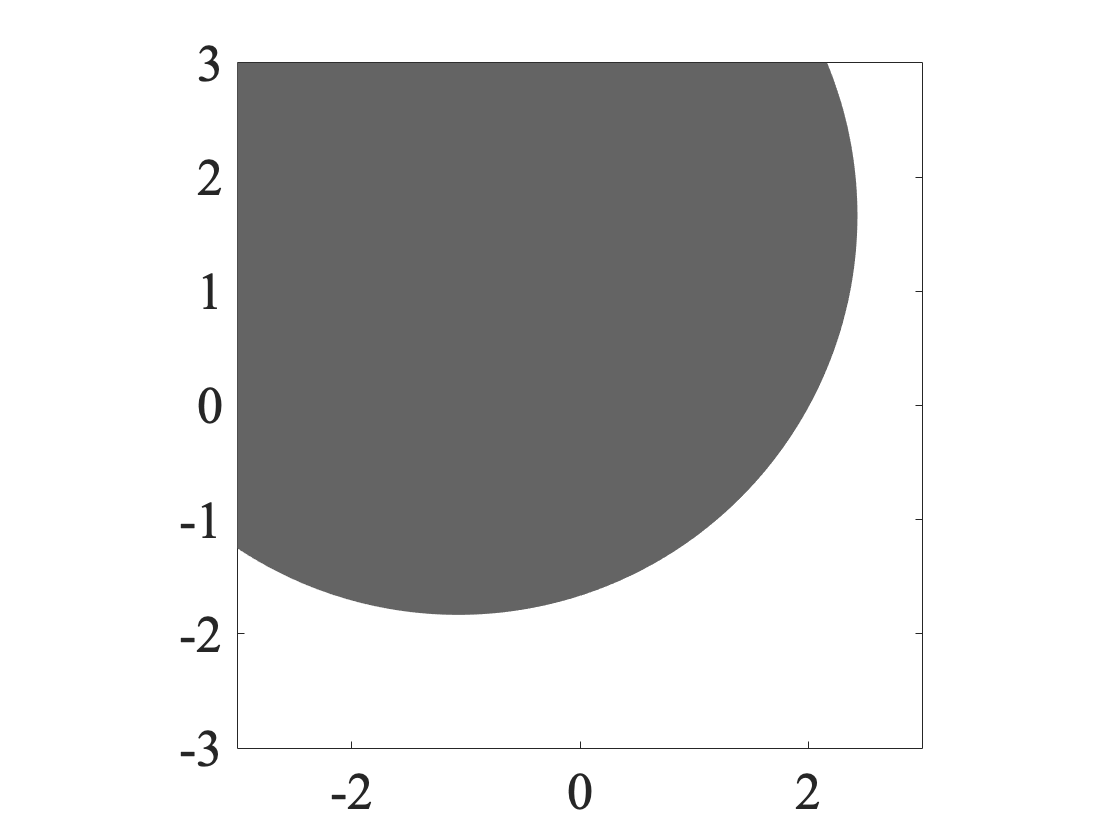}
		\caption{$\batSiz = 1$}
		\label{fig:batSizShapes12}
	\end{subfigure}
	\begin{subfigure}[b]{.19\linewidth}
		\centering \includegraphics[width=1\linewidth]{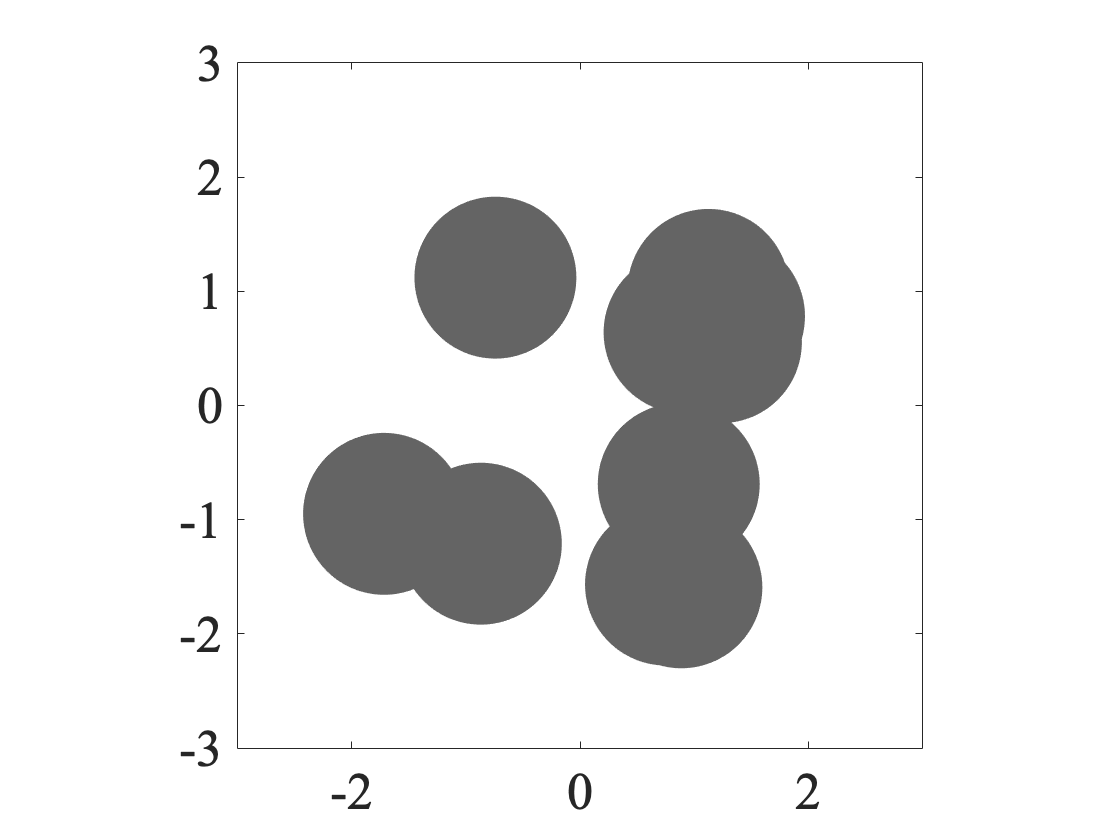}
		\caption{$\batSiz = 10$}
		\label{fig:batSizShapes102}
	\end{subfigure}
	\begin{subfigure}[b]{.19\linewidth}
		\centering \includegraphics[width=1\linewidth]{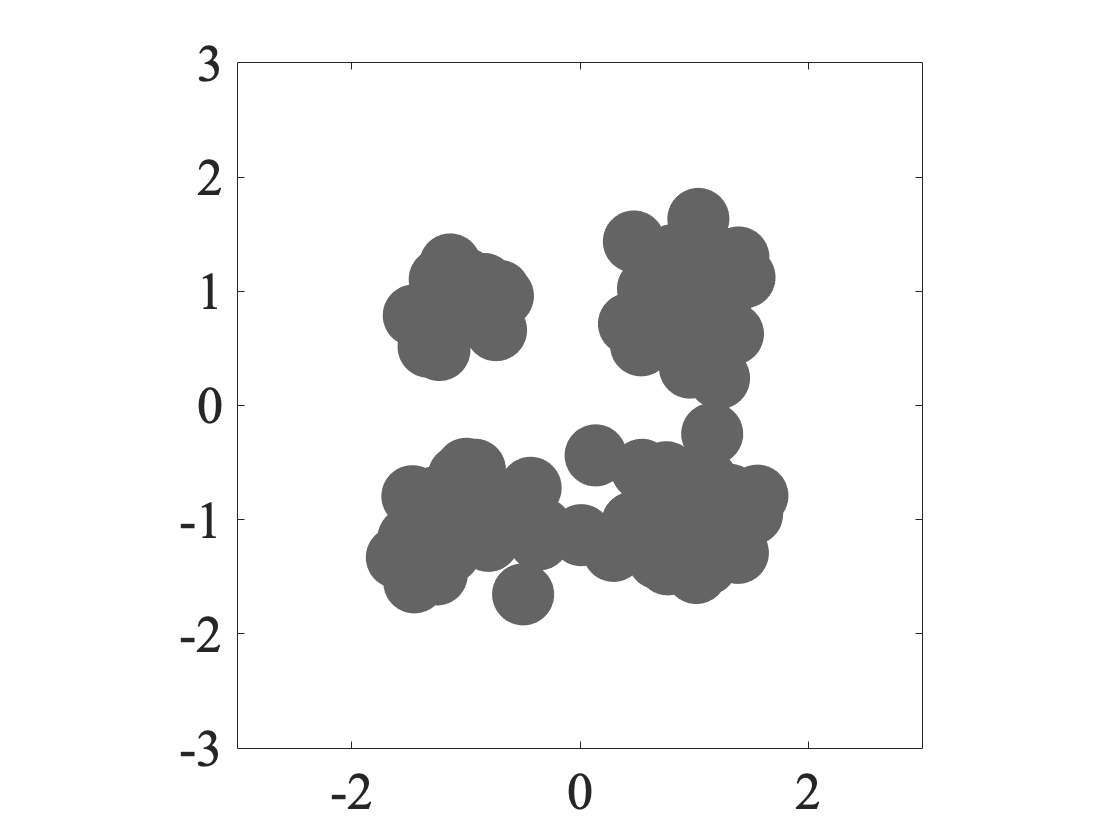}
		\caption{$\batSiz = 100$}
		\label{fig:batSizShapes1002}
	\end{subfigure}
	\begin{subfigure}[b]{.19\linewidth}
		\centering \includegraphics[width=1\linewidth]{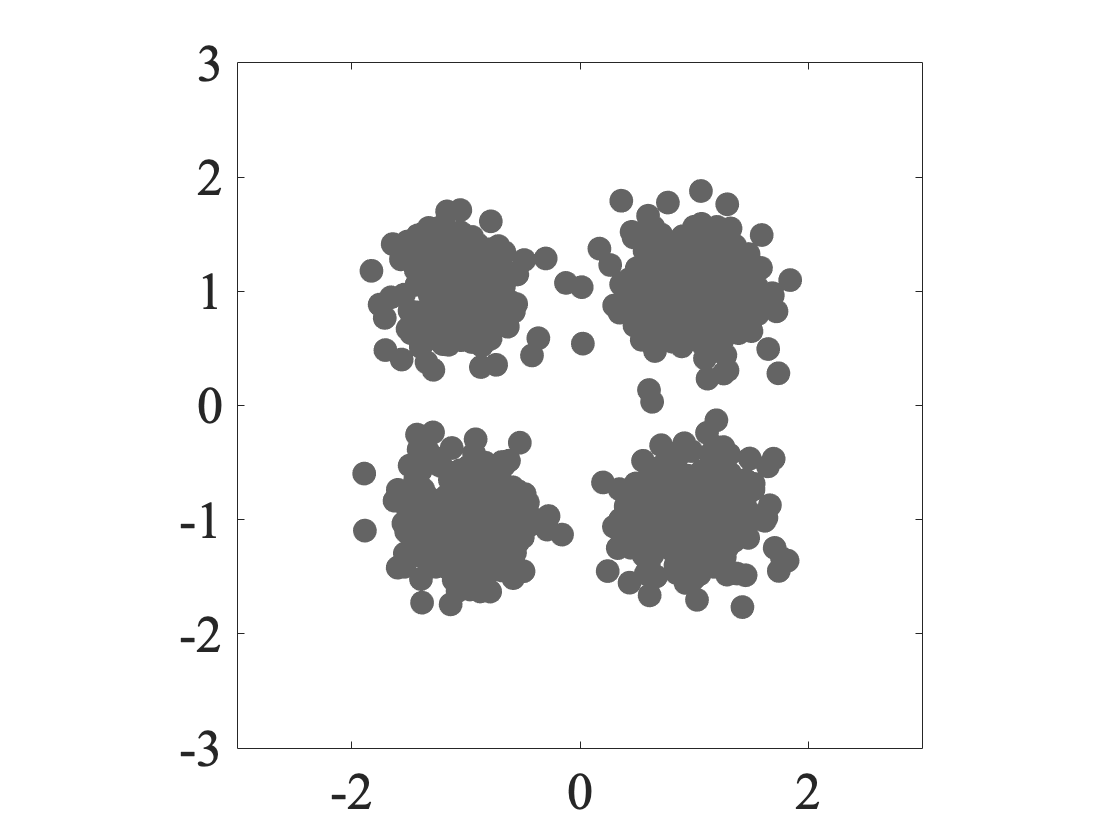}
		\caption{$\batSiz = 1000$}
		\label{fig:batSizShapes10002}
	\end{subfigure}

	\caption{The plots in the first column depict the probability density functions of three different $\Rset^2$-valued random vectors. The plots in columns two through five depict realizations of the uncertainty set $\lev{ F_n^{-1}(\conf_n)}$ for different values of $\batSiz$ for each random vector.}
	\label{fig:batSizShapes2}
\end{figure*}

\section{Experiments} \label{sec:experiments}

In this section, we present numerical experiments illustrating certain features and properties of the proposed uncertainty set estimator. 
Throughout these experiments, we consider three different random vectors, each of which is distributed according to a distinct Gaussian mixture. 

\subsection{Examining the Role of $m$}\label{sec:tradeoffs}

Figure \ref{fig:batSizShapes2} illustrates the impact of the number of balls $\batSiz$ on the shape and volume of the uncertainty sets produced by our method for three different  random vectors distributed according to Gaussian mixtures whose density functions are depicted in Figures \ref{fig:density3}, \ref{fig:density}, and \ref{fig:density2}.
In these experiments, we consider a fixed probability mass $\relLev = 0.9$, tolerance parameter $\buf = 0.05$, and confidence parameter $\confLev = 0.05$. In specifying the uncertainty set estimator, we set $\relLev_n = \relLev + \lambda^{*} \buf$ for all $n$, where  $\lambda^{*}$ is specified according to \eqref{eq:optimalLambda}. For each value of $m \in \{1,10,100, 1000\}$, we construct an uncertainty set $\lev{ F_n^{-1}(\conf_n)}$ according to \eqref{eqn:unc_set_def} using a training sample with size satisfying the bound \eqref{eq:samplesize}.

As expected, the plots in Figure \ref{fig:batSizShapes2}
show that the proposed quantile estimate for the $\alpha$-covering level, (i.e., the radius of each ball) appears to decrease with the number of balls used in the approximation of the $\alpha$-probability region.
And, on balance, the volume of the uncertainty sets appears to shrink with $m$, resulting in a potential decrease in the conservatism of the approximation to the original chance constrained problem. In particular, notice that the volume of the uncertainty set ($m=1$) depicted in Figure \ref{fig:batSizShapes12} is much larger than the volume of the uncertainty set ($m=1000$) depicted in Figure \ref{fig:batSizShapes10002}. From Figure \ref{fig:density2}, we see that the probability mass of the underlying distribution in this example is concentrated at four disconnected regions. Due to the separation between these high probability regions, uncertainty sets constructed according to a single ball are likely to have large volume. By contrast, notice that the difference in volume between the uncertainty sets in Figures \ref{fig:batSizShapes1} and \ref{fig:batSizShapes1000} is smaller, because in that example the underlying distribution has a significant amount of probability mass concentrated in a small connected region. Furthermore, the shapes and volumes of the uncertainty sets in Figures \ref{fig:batSizShapes13}-\ref{fig:batSizShapes10003}  do not vary significantly  for different values of $m$, because in this case the vast majority of the probability mass is concentrated in a circular region. It is important to note, however, that as $m$ increases, the volume of outliers also appears to increase, potentially resulting in increased conservatism of the approximation. Thus, from a practical perspective, the number of balls used in the approximation should be tuned to balance this tradeoff. 

\begin{figure}[htb!]
	\centering \includegraphics[width=.75\linewidth]{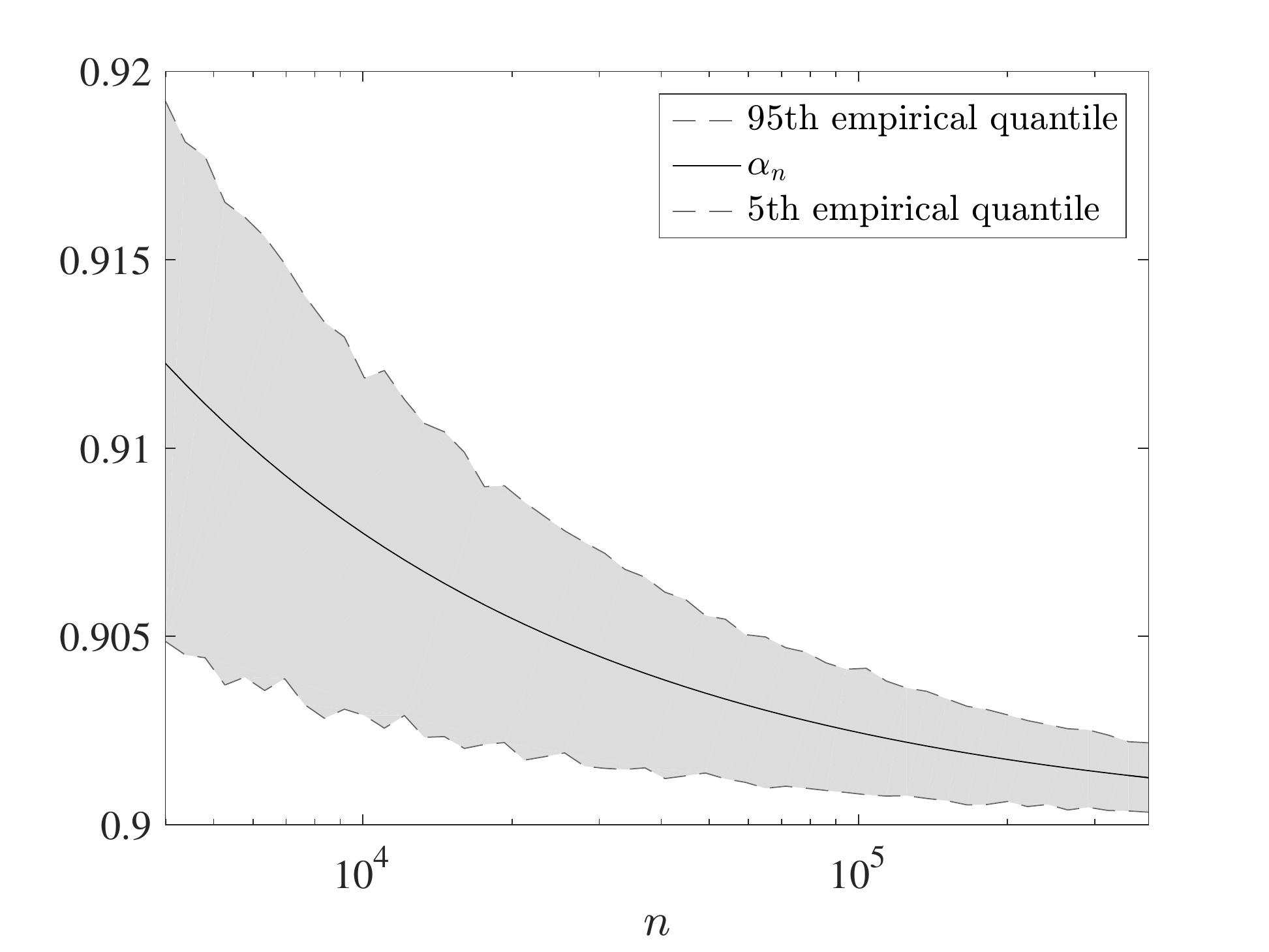}
	\caption{The middle 90\% empirical confidence interval associated with the coverage probability $\pi_\phi \left( F_n^{-1}(\conf_n) \right)$  for different tolerance parameter values.}
	\label{fig:difBuf}
\end{figure}

\subsection{Illustrating  Mass Consistency }\label{sec:consistency}
In Figure \ref{fig:difBuf}, we conduct a Monte-Carlo analysis to illustrate the mass consistency of the proposed uncertainty set estimator by varying the tolerance parameter value. Throughout these experiments, we fix the number of balls to $\batSiz = 10$, consider a target probability mass of $\relLev = 0.9$, and consider a random vector distributed according to a Gaussian mixture whose density is depicted in Figure \ref{fig:density}. We vary the tolerance parameter on a logarithmic scale between 0.005 and  0.05, and set $\relLev_\batNum =  \alpha+\lambda^{*} \buf$, where $\lambda^{*}$ is defined according to \eqref{eq:optimalLambda}. For each value of the tolerance parameter, we chose the number of training samples to be the smallest integer satisfying the sample size requirement \eqref{eq:samplesize}.

The numerical experiments are conducted as follows for each value of the tolerance parameter $\buf$. 
First, we generate a random shape sample $\bar{\uncPar}_{1}, \dots, \bar{\uncPar}_{\batSiz}$ that determines the centers of the balls composing the uncertainty set that we construct. 
For each experiment, we draw an i.i.d. training sample to evaluate $ F_n^{-1}(\conf_n)$ according to \eqref{eq:quantile} and $\lev{ F_n^{-1}(\conf_n)}$ according to \eqref{eqn:unc_set_def}. We then estimate the probability mass covered by each uncertainty set $\pi_\phi \left( F_n^{-1}(\conf_n) \right)$ using an empirical average based on one million independent samples of the random vector $\rand$. We estimate empirical confidence intervals associated with  $\pi_\phi \left( F_n^{-1}(\conf_n) \right)$ using one thousand independent experiments. 
Figure \ref{fig:difBuf} depicts the middle 90\% empirical confidence interval associated with the probability mass covered by the uncertainty sets for each value of $\buf$. Notice that, consistent with the result in Theorem \ref{thm:bound}, the fifth percentile of the probability mass covered by the uncertainty set remains above the target probability mass. Furthermore, the mass captured by the uncertainty set approaches the target probability mass with high confidence as the tolerance parameter value decreases to zero.

\section{Conclusion} \label{sec:conclusion}

 We provide a data-driven method to construct uncertainty sets for robust optimization problems. The estimators we consider are consistent, satisfy finite-sample performance guarantees, are efficient to compute, and give rise to tractable robust constraints for a large family of constraint functions. Furthermore, the proposed method provides the user with two mechanisms by which to control the complexity of the resulting robust optimization problem. First, the user can directly control the number of constraints in the resulting approximation while preserving the conservatism of the approximation with respect to the original chance constraint. Second, the geometry of the uncertainty sets can be tailored to accommodate the structure of the  given constraint function to ensure the computational tractability of the robust approximation. As a direction for future research, it would be interesting to refine our sample complexity results to explicitly reflect the role of the number of balls used in the approximation in determining the volume of the sets learned.

\bibliographystyle{IEEEtran}
\bibliography{references}{\markboth{References}{References}}

\end{document}